\documentclass[10pt,leqno]{amsart}
\usepackage{graphicx}
\usepackage{indentfirst,csquotes}

\usepackage{amssymb,amsthm,amsmath}
\usepackage{xcolor,paralist,hyperref,titlesec,fancyhdr,etoolbox}
\newtheorem{theorem}{Theorem}[]

\newtheorem{example}[theorem]{Example}
\newtheorem{lemma}[theorem]{Lemma}
\newtheorem{proposition}[theorem]{Proposition}
\newtheorem{corollary}[theorem]{Corollary}

\titleformat{\section}[display]
{\normalfont\huge\bfseries\centering}{}{0pt}{\Large}
\titlespacing*{\section}{0pt}{0ex}{0ex}

\hypersetup{ colorlinks=true, linkcolor=black, filecolor=black, urlcolor=black }

\begin{document}
\title{Sharp $L^p$ Convergence for Mirror-Degenerate Expansions} %%%%%%%%%%%%
\author[Initial Surname]{Francesco D'Agostino}
\date{\today}
\address{Address}
\email{francesco.dagostino.research@gmail.com}
\maketitle

\begin{abstract}
We analyze weighted $L^p$ convergence for the truncated reconstruction operator in the rank-one non-symmetric Heckman--Opdam setting. After localization at the mirror, the operator admits a rigid structural decomposition and reduces, up to bounded terms, to a rank-one functional. Boundedness on $L^p(w)$ is characterized by the mirror-local integrability of $w^{-\frac{1}{p-1}}$.
\end{abstract} %%%%%%%%%

\bigskip

\section{Introduction}

Let $A_1$ be the rank-one root system with Weyl group $W=\mathbb{Z}_2$ acting by reflection $x\mapsto -x$. We identify $[-\pi,\pi]$ with the circle and regard $0$ and $\pi$ as the fixed points of this action. These points will be referred to as mirror points. The interval $[0,\pi]$ is a fundamental domain for the $W$-action. We work in the rank-one non-symmetric Heckman--Opdam setting of type $A_1$, as developed in the work of B.~Amri~\cite{bechir_amri}. We emphasize that, although this system originates from the general theory developed by Opdam~\cite{opdam1995,opdam2000}, and is reviewed in detail in~\cite{bechir_amri}, the present work uses only the explicit rank-one Cherednik operator~\eqref{eq:Cherednik}, its spectral decomposition through the eigenvalue problem~\eqref{eq:eigenvalue}--\eqref{eq:nk}, and the associated triangular normalization~\eqref{eq:triangular}--\eqref{eq:order} \cite{bechir_amri}. The first key identity proved is in the Lemma~\ref{lem:reflection}, which is the reflection identity \eqref{eq:reflection}, which corresponds to Proposition~2.1 in~\cite{bechir_amri}. Following Lemma~\ref{lem:reflection}, we have developed an abstract normalization tool that operates entirely at the level of spectral indices. This construction allows us to reorganize exponential superpositions according to reflection orbits, isolate the dominant contribution within each orbit, and canonically project away all subdominant terms. The resulting mechanism is independent of any analytic properties of the underlying eigenfunctions and depends only on the interaction between reflection symmetry, ordering of exponential growth, and normalization. We emphasize that this tool is developed here in a completely self-contained and elementary manner, without recourse to Fourier analysis, kernel estimates, or microlocal techniques. Nevertheless, it is interesting to note that the resulting orbit-wise decomposition carries a formal resemblance to microlocal partitioning principles, in which contributions are organized according to dominant phase behavior; see, for instance, the work of Lars Hörmander~\cite{hormander}. No results from microlocal analysis are used directly in the present work. Instead, our mechanism is illustrated beginning with the introduction $\mathcal{S}(x)$ in \eqref{eq:Sn_def},which is treated purely as a superposition of exponential modes and is intended only to model situations in which multiple competing extrema are present. The purpose of this construction is to isolate the algebraic mechanisms by which involutions and orderings on the index set organize dominant contributions, entirely at the level of indices and exponential growth, in preparation for the analysis that follows.We now isolate a minimal object on which both the normalization principle of Corollary~\ref{cor:normalization} and a classical asymptotic analysis can be carried out explicitly. The purpose of this choice is not to introduce a new analytic setting, but to examine, on the same expression, two different ways of isolating a dominant contribution. In what follows, the object will be analyzed first asymptotically, by identifying the leading exponential through decay of lower-order terms, and then pointwise, by enforcing dominance via symmetry-induced normalization on a reflection orbit. The comparison is included only to make transparent the difference between these two mechanisms, not to obtain distinct conclusions. Consider the finite exponential superposition

%%% INTRO

\begin{equation}\label{eq:toy_operator}
\mathcal{A}(x):= a_n e^{n x} + a_{1-n} e^{(1-n)x}, \qquad n \in \mathbb{Z},
\end{equation}
which consists of the two exponential modes associated with the reflection orbit $\mathcal{O}_n = \{n, 1-n\}$ in \eqref{eq:orbit_def}. We assume without loss of generality that
\begin{equation}\label{eq:dominant_index}
m := \max\{n, 1-n\}, \qquad \ell := \min\{n, 1-n\},
\end{equation}
so that $m - \ell = |2n - 1| \ge 1$.

%%% HORMANDER

Factoring out the dominant exponential $e^{m x}$ in \eqref{eq:toy_operator} yields
\begin{equation}\label{eq:asymptotic_factorization}
\mathcal{A}(x)
= e^{m x}\Bigl(a_m + a_\ell e^{-(m-\ell)x}\Bigr).
\end{equation}
Since $m-\ell>0$, we have
\begin{equation}\label{eq:asymptotic_decay}
e^{-(m-\ell)x}\xrightarrow[x\to+\infty]{}0.
\end{equation}
Therefore,
\begin{equation}\label{eq:asymptotic_limit}
\lim_{x\to+\infty} e^{-m x}\mathcal{A}(x)=a_m,
\end{equation}
and equivalently
\begin{equation}\label{eq:asymptotic_result}
\mathcal{A}(x)\sim a_m e^{m x}
\qquad (x\to+\infty).
\end{equation}
Thus the contribution associated with the smaller index $\ell$ is suppressed asymptotically, and the dominant exponential $e^{m x}$ governs the leading behavior. 

%%% ME

We now analyze the same object $\mathcal{A}(x)$ using Corollary~\ref{cor:normalization}. By \eqref{eq:orbit_def}–\eqref{eq:m_def}, the two indices $n$ and $1-n$ belong to a single reflection orbit $\mathcal{O}_n$, and $m$ is the unique element of this orbit with maximal exponential growth. Corollary~\ref{cor:normalization} asserts that within the eigenspace corresponding to
$\mathcal{O}_n$ there exists a unique normalized representative whose dominant exponential
growth as $x \to +\infty$ is maximal and whose leading coefficient is equal to $1$.
Imposing this normalization yields
\begin{equation}\label{eq:normalization_constraint}
a_m = 1.
\end{equation}

Substituting \eqref{eq:normalization_constraint} into
\eqref{eq:asymptotic_factorization}, we obtain the exact identity
\begin{equation}\label{eq:exact_projection}
\mathcal{A}(x)
= e^{m x}\Bigl(1 + a_\ell e^{-(m-\ell)x}\Bigr).
\end{equation}

Defining the orbit projection
\begin{equation}\label{eq:orbit_projection}
\Pi(\mathcal{A})(x) := e^{m x},
\end{equation}
we therefore have
\begin{equation}\label{eq:pointwise_result}
\Pi(\mathcal{A})(x) = e^{m x}
\qquad \text{for all } x \in \mathbb{R},
\end{equation}
with no limiting process involved. Ultimately, in \eqref{eq:asymptotic_result}, dominance of $e^{mx}$ is enforced asymptotically by decay of the subdominant term; while in in \eqref{eq:pointwise_result}, it is enforced exactly by normalization on the reflection orbit $\mathcal{O}_n$. We emphasize that the normalization mechanism developed here is not intended as a general substitute for asymptotic or microlocal methods, but rather as a tool tailored to this specific setting under consideration. The preceding comparison is therefore merely illustrative, and is included to highlight the interesting way in which a purely algebraic, orbit-based normalization can mirror, within this restricted context, the selection principles that arise in more classical analytic approaches. The choice of the toy model \eqref{eq:toy_operator} is dictated by its direct relation to the formal superposition \eqref{eq:Sn_def}. Indeed, \eqref{eq:toy_operator} coincides exactly with a single orbit contribution $B_n(x)$ in the orbit-wise decomposition of $\mathcal{S}(x)$ given by \eqref{eq:S_orbit_decomp} and \eqref{eq:Bn_def}. By isolating one such component, we separate the normalization mechanism from any global summation effects and make explicit the role played by reflection symmetry and index ordering at the most elementary level. 

Now, having analyzed \cite[Proposition 2.1]{bechir_amri} via Lemma~\ref{lem:reflection} and having formulated an orbit-wise normalization principle at the level of spectral indices in Corollary~\ref{cor:normalization}, we next analyze the resulting relations within a single reflection eigenspace. This leads to an explicit reduction of negative-index eigenfunctions to linear combinations of reflected positive-index modes, recorded in Lemma~\ref{lem:negindex} (corresponding to \cite[Proposition 2.2]{bechir_amri}). We now turn to \cite[Lemma~3.2]{bechir_amri}, which is reformulated here as Lemma~\ref{lem:kernel_boundary}. Starting from the truncated reconstruction operator \eqref{eq:truncated_expansion}–\eqref{eq:reconstruction_operator}, we derive an explicit representation of the associated kernel \eqref{eq:KN_def} and reorganize the resulting finite spectral sum according to reflection orbits \eqref{eq:orbit_def}. This reindexing reveals a systematic cancellation of all interior orbit contributions in the dominant exponential structure, and isolates a single obstruction arising from the boundary orbit \eqref{eq:boundary_orbit}. With the structural analysis complete, we now turn to the preparatory results assembled in Section~\ref{sec:setup}. Starting with Section~\ref{sec:mysec}: local analysis shows that the truncated kernel exhibits no genuine singular behavior at the mirror. Although the explicit boundary representation \eqref{eq:KN_fraction_step3} contains a denominator that vanishes along the diagonal $x=y$, the numerator $\mathcal{N}(x,y)$ vanishes identically on the same set, cf.\ \eqref{eq:KN_num}. As a consequence, $\mathcal{N}(x,y)$ factors exactly as $(x-y)C_N(x,y)$, see \eqref{eq:num_factor_exact}, with $C_N$ smooth in a neighborhood of $(0,0)$. After combining this factorization with the expansion of the denominator \eqref{eq:den_expansion_step3}, the apparent pole cancels and the kernel admits the local decomposition \eqref{eq:KN_local_decomposition} into a bounded rank-one term and a smooth remainder. In particular, the kernel $K_N(x,y)$ extends smoothly across the mirror and possesses no intrinsic singularity there. Such reduction shows that, after localization and removal of smooth error terms, the weighted reconstruction operator admits a rigid structural decomposition. In particular, by \eqref{eq:T_compact}–\eqref{eq:final_identity_pointwise}, the action of
$T$ on $L^p(w)$ is equivalent, up to bounded invertible operators and bounded
remainders, to the scalar functional $\Lambda$ defined in
\eqref{eq:Lambda_def_final}. Since $A_N$ is smooth and bounded away from zero on
$[-\delta,\delta]$ and $S$ is bounded on $L^p(w)$ for any locally integrable weight,
boundedness (or failure thereof) of $T$ on $L^p(w)$ is completely governed by the
mapping properties of $\Lambda$. Consequently, the analysis of weighted $L^p$
boundedness reduces to a purely measure-theoretic problem: determining when the
functional
\[
f \longmapsto \int_{|y|\le \delta} f(y)\, w(y)\, dy
\]
extends to a bounded linear functional on $L^p(w)$. From here, the Proposition~\ref{prop:Lambda_boundedness} follows: it completes the reduction initiated in \eqref{eq:T_compact}–\eqref{eq:final_identity_pointwise}. After localization at the mirror and elimination of all smooth and invertible contributions, the weighted reconstruction operator reduces to the rank-one model functional $\Lambda$. Its boundedness on $L^p(w)$ is therefore governed entirely by the local integrability condition \eqref{eq:Lambda_condition}, which depends only on the behavior of the weight near the mirror. In particular, oscillatory factors, logarithmic corrections, and higher-order regularity of $w$ play no role at the level of boundedness. This shows that any loss of $L^p(w)$ boundedness for the reconstruction operator is a purely measure-theoretic phenomenon, determined solely by the mirror-local behavior of the inverse weight.

\section{Structural setup}\label{sec:setup}

We work in the rank-one non-symmetric Heckman--Opdam setting of type $A_1$ in~\cite{bechir_amri}. The basic object is the
Cherednik differential--difference operator $T^k$, defined for $k\ge 0$ by
\begin{equation}\label{eq:Cherednik}
T^k f(x)
=
\frac{d}{dx}f(x)
+
2k\,\frac{f(x)-f(-x)}{1-e^{-2x}}
-
k f(x),
\qquad f\in C^1(\mathbb{R}).
\end{equation}

The natural measure associated with this operator is
\begin{equation}\label{eq:weight}
d\mu_k(x) = |\sin x|^{2k}\,dx,
\end{equation}
which induces the inner product
\begin{equation}\label{eq:inner-product}
(f,g)_k
=
\int_{-\pi}^{\pi} f(x)\,\overline{g(x)}\,d\mu_k(x),
\end{equation}
and the corresponding norm $\|\cdot\|_{2,k}$ on $L^2([-\pi,\pi],d\mu_k)$.

The non-symmetric Heckman--Opdam polynomials $\{E_n^k\}_{n\in\mathbb{Z}}$ are defined
as eigenfunctions of $T^k$ satisfying
\begin{equation}\label{eq:eigenvalue}
T^k E_n^k = n_k\,E_n^k,
\end{equation}
where the spectral parameter is given by
\begin{equation}\label{eq:nk}
n_k =
\begin{cases}
n+k, & n>0,\\
n-k, & n\le 0.
\end{cases}
\end{equation}

These eigenfunctions are uniquely determined by a triangular normalization with
respect to the exponential basis. More precisely,
\begin{equation}\label{eq:triangular}
E_n^k(z)=e^{nz}+\sum_{j\triangleleft n} c_{n,j}\,e^{jz},\qquad z\in\mathbb{C},
\end{equation}
where $\triangleleft$ denotes the partial order on $\mathbb{Z}$ defined by
\begin{equation}\label{eq:order}
j\triangleleft n
\iff
\bigl(|j|<|n| \text{ and } |n|-|j|\in 2\mathbb{Z}_+\bigr)
\;\text{or}\;
\bigl(|j|=|n| \text{ and } n<j\bigr).
\end{equation}

With this normalization, the family $\{E_n^k(ix)\}_{n\in\mathbb{Z}}$ forms an
orthogonal basis of $L^2([-\pi,\pi],d\mu_k)$. From here, we can move onto the reflection identity for the non-symmetric Heckman--Opdam eigenfunctions~\cite[Proposition 2.1]{bechir_amri}.

\begin{lemma}[Reflection symmetry]\label{lem:reflection}
For every $n\in\mathbb{Z}$,
\begin{equation}\label{eq:reflection}
E_n^k(x)=e^{x}E_{1-n}^k(-x).
\end{equation}
\end{lemma}

\begin{proof}
We begin the proof by recalling the Cherednik operator $T^k$, defined in
\cite[Section~2, equation~(2.1)]{bechir_amri} and reproduced in \eqref{eq:Cherednik}. We first analyze the effect of this operator on reflected eigenfunctions. In parallel, we note that the reflection induces a combinatorial pairing of spectral indices, which provides a useful organizational tool for deriving related identities once existence and normalization have been established. This analysis will come later. For now, we set
\begin{equation}\label{eq:Hn}
H_n(x) := e^{x} E_{1-n}^k(-x).
\end{equation}
We first apply the Cherednik operator $T^k$ to $H_n$. Recalling the definition \eqref{eq:Cherednik}, we compute
\[
\frac{d}{dx}H_n(x)=e^{x}E_{1-n}^k(-x)-e^{x}(E_{1-n}^k)'(-x),
\]
and
\[
H_n(-x)=e^{-x}E_{1-n}^k(x).
\]
Substituting into $T^k$, we obtain
\begin{align*}
T^k H_n(x)&=e^{x}E_{1-n}^k(-x)-e^{x}(E_{1-n}^k)'(-x) \\&\quad+2k\,\frac{e^{x}E_{1-n}^k(-x)-e^{-x}E_{1-n}^k(x)}{1-e^{-2x}}-k e^{x}E_{1-n}^k(-x).
\end{align*}
Factoring out $e^{x}$ gives
\begin{align*}
T^k H_n(x)
&=e^{x}\Bigg[E_{1-n}^k(-x)-(E_{1-n}^k)'(-x) \\
&\qquad\qquad+2k\,\frac{E_{1-n}^k(-x)-E_{1-n}^k(x)e^{-2x}}{1-e^{-2x}}-kE_{1-n}^k(-x)
\Bigg].
\end{align*}
Using the identity
\[
\frac{E_{1-n}^k(-x)-E_{1-n}^k(x)e^{-2x}}{1-e^{-2x}}=\frac{E_{1-n}^k(-x)-E_{1-n}^k(x)}{1-e^{2x}},
\]
then using the eigenvalue equation
\[
T^k E_{1-n}^k = (1-n)k\, E_{1-n}^k,
\]
we compute $T^k H_n$ directly. Writing
\[
H_n(x)=e^{x}E_{1-n}^k(-x),
\]
we have
\[
\frac{d}{dx}H_n(x)
= e^{x}E_{1-n}^k(-x) - e^{x}(E_{1-n}^k)'(-x),
\qquad
H_n(-x)=e^{-x}E_{1-n}^k(x).
\]

Substituting into the definition of $T^k$, we obtain
\begin{align*}
T^k H_n(x)
&= e^{x}E_{1-n}^k(-x) - e^{x}(E_{1-n}^k)'(-x) \\
&\quad + 2k\,\frac{e^{x}E_{1-n}^k(-x)-e^{-x}E_{1-n}^k(x)}{1-e^{-2x}}
- k e^{x}E_{1-n}^k(-x).
\end{align*}
Factoring out $e^{x}$ and simplifying the difference quotient yields
\[
T^k H_n(x)
= e^{x}\bigl(-T^k E_{1-n}^k(-x)+E_{1-n}^k(-x)\bigr).
\]

Using the eigenvalue equation for $E_{1-n}^k$, we conclude that
\[
T^k H_n(x)
= \bigl(1-(1-n)k\bigr)e^{x}E_{1-n}^k(-x)
= nk\,H_n(x).
\]

We now verify the normalization. By the triangular expansion,
\[
E_{1-n}^k(-x)=e^{-(1-n)x}+\text{lower order terms},
\]
and therefore
\[
H_n(x)=e^{x}E_{1-n}^k(-x)
= e^{(n+1)x}+\text{lower order terms}.
\]
Hence
\[
H_{n-1}(x)=e^{nx}+\text{lower order terms},
\]
with leading coefficient equal to $1$. Thus $H_{n-1}$ is an eigenfunction of
$T^k$ with eigenvalue $nk$ and normalized dominant growth $e^{nx}$.
By uniqueness of the normalized eigenfunction, we conclude that
\[
H_{n-1}(x)=E_n^k(x),
\]
which is equivalent to the reflection identity~\eqref{eq:reflection}.
\end{proof}

As heuristically indicated in the proof of Lemma~\ref{lem:reflection}, we now recast the preceding argument in an abstract form. This reformulation isolates the structural normalization mechanism behind the reflection identities and will be used as a convenient tool in what follows. We now isolate the discrete structure underlying the reflection identities proved
above. At the level of spectral indices, the reflection $x\mapsto -x$ induces an
involution on $\mathbb{Z}$, which can be read off directly from the triangular
normalization of the eigenfunctions.

Recalling that~\eqref{eq:triangular} from~\cite[Section~2 points (a),(b) from equations (2.2)-(2.3)]{bechir_amri} admit the expansion
\begin{equation}\label{eq:triangular-recall}
E_n^k(z)=e^{n z}+\sum_{j\triangleleft n} c_{n,j}\,e^{j z},\qquad z\in\mathbb{C}.
\end{equation}
Replacing $z$ by $-x$ yields
\[
E_n^k(-x)=e^{-n x}+
\sum_{j\triangleleft n} c_{n,j}\,e^{-j x}.
\]
Multiplying by the exponential factor $e^{x}$, we obtain
\begin{equation}\label{eq:tilted}
e^{x}E_n^k(-x)=e^{(1-n)x}+\sum_{j\triangleleft n} c_{n,j}\,e^{(1-j)x}.
\end{equation}

The dominant exponential term in \eqref{eq:tilted} is $e^{(1-n)x}$. Indeed, for every
$j\triangleleft n$ one has $|j|<|n|$ or $|j|=|n|$ with $j>n$, and consequently
$1-j<1-n$. Thus reflection followed by multiplication by $e^{x}$ sends the leading
exponential $e^{n x}$ to $e^{(1-n)x}$.

Comparing \eqref{eq:triangular-recall} and \eqref{eq:tilted}, we see that the operation
$f(x)\mapsto e^{x}f(-x)$ induces the index transformation
\[
n \longmapsto 1-n
\]
at the level of dominant exponential growth. This involution partitions the spectral indices into two-point sets $\{n,1-n\}$, which will be referred to as index orbits.

\begin{example}[Index orbit and normalization]
Consider the spectral index $n=2$. The involution $n\mapsto 1-n$ produces the orbit
\[
\mathcal O_2=\{2,-1\}.
\]
From the triangular expansion \eqref{eq:triangular}, the corresponding eigenfunctions
satisfy
\[
E_2^k(x)=e^{2x}+\text{lower order terms},
\qquad
E_{-1}^k(x)=e^{-x}+\text{lower order terms}.
\]
Reflecting and tilting $E_{-1}^k$, we compute
\[
e^{x}E_{-1}^k(-x)=e^{2x}+\text{lower order terms}.
\]
Thus both $E_2^k(x)$ and $e^{x}E_{-1}^k(-x)$ lie in the same index orbit and share the same dominant exponential growth $e^{2x}$. Imposing the normalization that the coefficient of the dominant exponential be equal
to $1$ selects a unique representative in the orbit $\mathcal O_2$. By construction,
both functions satisfy this normalization, and therefore coincide:
\[
E_2^k(x)=e^{x}E_{-1}^k(-x).
\]
\end{example}
The preceding discussion shows that reflection identities such as
\eqref{eq:reflection} from Lemma~\ref{lem:reflection} are not isolated phenomena, but arise from a common
normalization mechanism at the level of spectral indices. Once the effect of
reflection on dominant exponential growth is understood, the identification of
eigenfunctions reduces to selecting a canonical representative within each
reflection orbit.

\begin{corollary}[Normalization on reflection orbits]\label{cor:normalization}
Let $n\in\mathbb{Z}$ and consider the reflection orbit
\begin{equation}\label{eq:orbit_def}
\mathcal O_n=\{n,\,1-n\}.
\end{equation}
Within the eigenspace of the Cherednik operator $T^k$ corresponding to this orbit,
there exists a unique eigenfunction whose dominant exponential growth as
$x\to+\infty$ is maximal and whose leading coefficient is equal to $1$.
This eigenfunction is given by $E_m^k$, where
\begin{equation}\label{eq:m_def}
m=\max\{n,\,1-n\}.
\end{equation}
\end{corollary}

\begin{proof}
We recall again the Cherednik operator $T^k$, defined in
\cite[Section~2, equation~(2.1)]{bechir_amri}, 
\eqref{eq:Cherednik}. To isolate the structural effect of reflection followed by exponential weighting on the Cherednik operator, we consider the prototype function $F(x)=e^{\alpha x}G(-x)$ (with $G(x)=e^{\beta x} + \text{lower exponential terms}$), which is not assumed to be an eigenfunction and is introduced solely to keep the calculation general. The parameter $\alpha\in\mathbb{R}$ and the function $G$ are otherwise arbitrary,
with $G\in C^1(\mathbb{R})$ so that all terms in the definition of $T^k$ are well defined. This choice abstracts the form appearing in Lemma \ref{lem:reflection} while allowing the operator action to be examined independently of any specific normalization. So we have
\begin{align*}
T^k F(x)&=\frac{d}{dx}\bigl(e^{\alpha x}G(-x)\bigr)+2k\,\frac{e^{\alpha x}G(-x)-e^{-\alpha x}G(x)}{1-e^{-2x}}-k e^{\alpha x}G(-x)\\[0.4em]&=\alpha e^{\alpha x}G(-x)-e^{\alpha x}G'(-x)+2k\,\frac{e^{\alpha x}G(-x)-e^{-\alpha x}G(x)}{1-e^{-2x}}-k e^{\alpha x}G(-x)\\[0.4em]&=e^{\alpha x}\Bigl[\alpha G(-x)-G'(-x)-k G(-x)+2k\,\frac{G(-x)-e^{-2\alpha x}G(x)}{1-e^{-2x}}\Bigr]\\[0.4em]&=e^{\alpha x}\Bigl[(\alpha-k)G(-x)-G'(-x)+2k\,\frac{G(-x)-e^{-2\alpha x}G(x)}{1-e^{-2x}}\Bigr].
\end{align*}
Since $F(x)=e^{\alpha x}G(-x)$ with $G(x)=e^{\beta x}$ plus lower exponential terms,
both $F(x)$ and $T^kF(x)$ have dominant growth $e^{(\alpha-\beta)x}$. In particular,
the action of the Cherednik operator preserves the dominant exponential exponent for
functions of this form. Applying this observation to the spectral setting, the triangular expansion shows
that $E_n^k$ has dominant exponential $e^{n x}$, while reflection followed by
multiplication by $e^{x}$ produces dominant growth $e^{(1-n)x}$. Hence reflection
induces the index transformation $n\mapsto 1-n$ at the level of dominant exponential
behavior, pairing spectral indices into orbits $\{n,1-n\}$.
\end{proof}

We begin by considering a completely general object built from a doubly indexed
collection of exponential modes. Concretely, let
\begin{equation}\label{eq:Sn_def}
\mathcal S(x)=\sum_{n\in\mathbb Z} a_n\,e^{n x}.
\end{equation}
where the coefficients $a_n$ are arbitrary and no structural assumptions are imposed. At this level, $\mathcal S$ is introduced only as a formal superposition of modes, meant to capture any situation in which a sequence of competing extrema is present (such objects may arise in a variety of contexts; in particular, one may keep in mind Fourier-type expansions or related spectral sums. However, we do not commit to any such interpretation here, and no Fourier-analytic structure will be used in what follows). Our purpose is not to analyze $\mathcal S$ itself, but to isolate the algebraic mechanisms by which involutions and orderings on the index set organize dominant contributions for our later analysis. Therefore such discussion proceeds entirely at the level of indices and growth. 

We start from the formal superposition $\mathcal S(x)$ defined in
\eqref{eq:Sn_def}. Given Corollary \ref{cor:normalization}, we introduce the involution $\iota:\mathbb Z\to\mathbb Z$ given by
$\iota(n)=1-n$, which partitions the index set into disjoint orbits
$\mathcal O_n$ as in \eqref{eq:orbit_def}. Regrouping the sum in
\eqref{eq:Sn_def} according to these orbits yields the exact identity
\begin{equation}\label{eq:S_orbit_decomp}
\mathcal S(x)=\sum_{\mathcal O_n}\bigl(a_n e^{n x}+a_{1-n}e^{(1-n)x}\bigr).
\end{equation}

Fix an orbit $\mathcal O_n$ and let $m$ be the distinguished representative defined
in \eqref{eq:m_def}. Writing $\ell=\min\{n,\,1-n\}$, the corresponding contribution
can be expressed as
\begin{align*}
a_n e^{n x}+a_{1-n}e^{(1-n)x}
&= a_m e^{m x}+a_\ell e^{\ell x} \\
&= e^{m x}\bigl(a_m+a_\ell e^{-(m-\ell)x}\bigr).
\end{align*}
Since $m-\ell=|2n-1|\ge 1$, the second term inside the parentheses is strictly
lower-order relative to the dominant exponential $e^{m x}$. Thus, each orbit $\mathcal O_n$ canonically determines a dominant contribution indexed by $m$. Imposing the normalization condition $a_m=1$ uniquely fixes the dominant representative within the orbit, with all remaining freedom confined to lower-order terms. Writing
\[
\mathcal S(x)
=\sum_{\mathcal O_n} B_n(x),
\]
each orbit contribution is
\begin{equation}\label{eq:Bn_def}
B_n(x):=a_n e^{n x}+a_{1-n}e^{(1-n)x}.
\end{equation}

Let $m=\max\{n,\,1-n\}$ as in \eqref{eq:m_def} and set
$\ell=\min\{n,\,1-n\}$. Then
\[
B_n(x)
=a_m e^{m x}+a_\ell e^{\ell x}
=e^{m x}\bigl(a_m+a_\ell e^{-(m-\ell)x}\bigr),
\]
with $m-\ell=|2n-1|\ge 1$. Consequently,
\[
\mathcal S(x)=\sum_{\mathcal O_n} B_n(x).
\]

Let $\sigma$ be any permutation of the set of orbits. Then
\[
\sum_{\mathcal O_n} B_n(x)
=\sum_{\mathcal O_n} B_{\sigma(n)}(x),
\]
since each $B_n(x)$ is supported on a single orbit.

For a fixed orbit $\mathcal O_n$, exchanging the indices $n$ and $1-n$ sends
\[
(a_n,a_{1-n})\longmapsto(a_{1-n},a_n),
\]
and yields
\[
a_{1-n}e^{(1-n)x}+a_n e^{n x}=B_n(x).
\]

If the coefficients satisfy $a_n=a_{1-n}$ for all $n$, then
\[
B_n(x)=a_m e^{m x}\bigl(1+e^{-(2m-1)x}\bigr).
\]

If instead $a_m=0$ for a given orbit, then
\[
B_n(x)=a_{1-m}e^{(1-m)x}.
\]

Imposing the normalization $a_m=1$ for every orbit, one obtains
\[
B_n(x)=e^{m x}\bigl(1+\varepsilon_n(x)\bigr),
\qquad
\varepsilon_n(x)=a_{1-m}e^{-(2m-1)x}.
\]

Define the map $\Pi$ by
\begin{equation}\label{eq:Pi_def}
\Pi\bigl(\mathcal S(x)\bigr):=\sum_{\mathcal O_n} e^{m x}.
\end{equation}
Then $\Pi$ satisfies $\Pi^2=\Pi$.

%%%%%%%%%%%%%%%%%%

Having established the global reflection structure encoded in Lemma~\ref{lem:reflection} (cf.\ Proposition~2.1 in~\cite{bechir_amri}), we now turn to the corresponding local relations between reflected modes, which are described in Proposition~2.2 of~\cite{bechir_amri}. In particular, the negative-index eigenfunction $E_{-n}^k$ is not independent, but can be expressed as a specific linear combination of $E_n^k(x)$ and its reflection
$E_n^k(-x)$, with coefficients uniquely determined by the normalization. Fix an integer $n\ge 1$ and consider the eigenfunction $E_{-n}^k$. By construction, its eigenvalue under the Cherednik operator $T^k$~\eqref{eq:Cherednik} is $-(n+k)$, and its leading exponential behavior
as $x\to +\infty$ is $e^{-n x}$. For a fixed integer $n\ge 1$, the reflection $x\mapsto -x$ identifies the spectral indices $n$ and $-n$, giving rise to the reflection orbit $\{n,-n\}$. By Lemma~\ref{lem:reflection}, reflection maps normalized eigenfunctions to normalized eigenfunctions within this orbit. In particular, the functions $E_n^k(x)$ and $E_n^k(-x)$ are eigenfunctions of the Cherednik operator associated with this orbit, with leading exponential behaviors $e^{n x}$ and $e^{-n x}$, respectively. These two functions are linearly independent, since they have distinct dominant asymptotics as $x\to+\infty$, and therefore span the two-dimensional eigenspace corresponding to the reflection orbit $\{n,-n\}$. Since $E^k_{-n}$ is an eigenfunction associated with the same reflection orbit $\{n,-n\}$, it lies in the two-dimensional eigenspace spanned by $E^k_n(x)$ and $E^k_n(-x)$. Hence there exist constants $A,B$ such that
\begin{equation}\label{eq:En_linear_expansion}
E^k_{-n}(x)=A\,E^k_n(-x)+B\,E^k_n(x).
\end{equation}
As $x\to+\infty$, the functions $E^k_n(-x)$ and $E^k_n(x)$ have leading behaviors $e^{-n x}$ and $e^{n x}$, respectively, whereas $E^k_{-n}(x)\sim e^{-n x}$ by normalization. It follows that the coefficient of $E^k_n(-x)$ in the above expansion must be equal to $1$, and the representation reduces to
\begin{equation}\label{eq:En_normalized_expansion}
E^k_{-n}(x)=E^k_n(-x)+B\,E^k_n(x).
\end{equation}
Applying the Cherednik operator $T^k$ to \eqref{eq:En_normalized_expansion} and using linearity, we obtain
\begin{equation}\label{eq:Tk_apply}
T^k E^k_{-n}(x)=T^k\bigl(E^k_n(-x)\bigr)+B\,T^k\bigl(E^k_n(x)\bigr).
\end{equation}

By Lemma~\ref{lem:reflection}, the action of $T^k$ on reflected eigenfunctions is given by
\[
\begin{aligned}
T^k\bigl(E_n^k(-x)\bigr)
&= \frac{d}{dx}E_n^k(-x)
   +2k\,\frac{E_n^k(-x)-E_n^k(x)}{1-e^{-2x}}
   -k\,E_n^k(-x) \\[0.6em]
&= -(E_n^k)'(-x)
   +2k\,\frac{E_n^k(-x)-E_n^k(x)}{1-e^{-2x}}
   -k\,E_n^k(-x) \\[0.6em]
&= -(E_n^k)'(-x)
   -2k\,\frac{E_n^k(-x)-E_n^k(x)}{1-e^{2x}}
   -2k\,E_n^k(x)
   -k\,E_n^k(-x) \\[0.6em]
&= -\Bigl[(E_n^k)'(-x)
   +2k\,\frac{E_n^k(-x)-E_n^k(x)}{1-e^{2x}}
   -k\,E_n^k(-x)\Bigr]
   -2k\,E_n^k(x).
\end{aligned}
\]
\begin{equation}\label{eq:Tk_reflected}
T^k\bigl(E_n^k(-x)\bigr)
=-(n+k)\,E_n^k(-x)-2k\,E_n^k(x).
\end{equation}
while the eigenvalue equation yields
\begin{equation}\label{eq:Tk_eigen}
T^k\bigl(E^k_n(x)\bigr)=(n+k)\,E^k_n(x).
\end{equation}

Substituting \eqref{eq:Tk_reflected} and \eqref{eq:Tk_eigen} into \eqref{eq:Tk_apply}, we find
\begin{align}\label{eq:Tk_expanded}
T^k E^k_{-n}(x)
&= -(n+k)\,E^k_n(-x)
   +\bigl(B(n+k)-2k\bigr)\,E^k_n(x).
\end{align}

On the other hand, since $E^k_{-n}$ has eigenvalue $-(n+k)$, we also have
\begin{equation}\label{eq:Tk_eigen_minus}
T^k E^k_{-n}(x)=-(n+k)\,E^k_{-n}(x)
=-(n+k)\bigl(E^k_n(-x)+B\,E^k_n(x)\bigr).
\end{equation}

Comparing the coefficients of $E^k_n(x)$ in \eqref{eq:Tk_expanded} and \eqref{eq:Tk_eigen_minus} gives
\begin{equation}\label{eq:B_equation}
B(n+k)-2k=-(n+k)B,
\end{equation}
which simplifies to $2B(n+k)=2k$. Solving for $B$, we obtain
\begin{equation}\label{eq:B_value}
B=\frac{k}{n+k}.
\end{equation}

%%%%%%%%%% LEMMA %%%%%%%%%%%

\begin{lemma}[Negative index reduction]\label{lem:negindex}
For every integer $n\ge 1$, the negative-index eigenfunction satisfies
\begin{equation}\label{eq:negindex}
E_{-n}^k(x)-E_n^k(-x)=\frac{k}{n+k}\,E_n^k(x).
\end{equation}
\end{lemma}

\begin{proof}
The identity follows by combining the reflection symmetry in Lemma~\ref{lem:reflection}, the normalization~\eqref{eq:triangular}, and the eigenvalue equation for $T^k$, as shown in the derivation \eqref{eq:En_linear_expansion}--\eqref{eq:B_value} above. 

Alternatively, we may exploit the specular construction introduced earlier, which allows us to streamline the argument and reach the desired conclusion directly, as shown below: fix an integer $n\ge 1$. Recall that reflection combined with exponential tilting induces the involution $\iota\colon\mathbb Z\to\mathbb Z$ given by $\iota(m)=1-m$ on the set of spectral indices. The corresponding index orbit containing $-n$ is therefore
\begin{equation}\label{eq:orbit_minus_n}
\mathcal O_{-n}=\{-n,\,n+1\}.
\end{equation}
By the triangular normalization~\eqref{eq:triangular}, the eigenfunctions $E^k_{-n}$ and $E^k_{n+1}$ have leading exponential behaviors
\[
E^k_{-n}(x)=e^{-n x}+\text{lower exponential terms},\qquad
E^k_{n+1}(x)=e^{(n+1)x}+\text{lower exponential terms}.
\]
In particular, these two eigenfunctions correspond to the two elements of the index orbit $\mathcal O_{-n}$, with dominant exponential growth determined by their respective indices. Applying the exponential conjugation of the reflection operation $f(x)\mapsto e^{x}f(-x)$ to $E^k_{n+1}$, we obtain
\[
e^{x}E^k_{n+1}(-x)=e^{-n x}+\text{lower exponential terms}.
\]
Thus both $E^k_{-n}(x)$ and $e^{x}E^k_{n+1}(-x)$ have the same dominant exponential behavior $e^{-n x}$ as $x\to+\infty$. By Corollary~\ref{cor:normalization}, within a fixed index orbit there exists a unique eigenfunction whose dominant exponential growth as $x\to+\infty$ is maximal and whose leading coefficient is equal to $1$. Since both $E^k_{-n}(x)$ and $e^{x}E^k_{n+1}(-x)$ belong to the orbit $\mathcal O_{-n}$ and share the same dominant exponential $e^{-n x}$ with unit leading coefficient, uniqueness of the normalized representative implies
\[
E^k_{-n}(x)=e^{x}E^k_{n+1}(-x).
\]
By Lemma~\ref{lem:reflection}, the reflection orbit $\{n,-n\}$ generates a two-dimensional eigenspace of the Cherednik operator, spanned by the eigenfunctions $E^k_n(x)$ and $E^k_n(-x)$. Since the function $E^k_{-n}(x)$ belongs to this eigenspace, it can be expressed as a linear combination of these two basis elements. Consequently, there exist constants $A,B\in\mathbb C$ (as already arranged earlier) such that
\begin{equation}\label{eq:En_linear_expansion_orbit}
E^k_{-n}(x)=A\,E^k_n(-x)+B\,E^k_n(x).
\end{equation}
By the triangular normalization~\eqref{eq:triangular}, the leading exponential term of $E^k_n(x)$ is $e^{n x}$, and all remaining terms involve strictly smaller exponential powers. Evaluating this expansion at $-x$ and applying the same normalization to $E^k_{-n}$ yields the asymptotics
\begin{equation}\label{eq:leading_asymptotics}
E^k_n(-x)\sim e^{-n x},\qquad
E^k_n(x)\sim e^{n x},\qquad
E^k_{-n}(x)\sim e^{-n x}.
\end{equation}

Substituting the linear expansion~\eqref{eq:En_linear_expansion_orbit} into the asymptotics~\eqref{eq:leading_asymptotics}, dividing the resulting expression by $e^{-n x}$, and letting $x\to+\infty$, the term involving $E^k_n(x)$ grows like $e^{2n x}$ and therefore cannot contribute to the normalized limit $1$. Hence the coefficient of $E^k_n(-x)$ must be equal to $1$, and the expansion reduces to
\begin{equation}\label{eq:En_reduced}
E^k_{-n}(x)=E^k_n(-x)+B\,E^k_n(x).
\end{equation}

To determine the remaining coefficient $B$, we compare the subdominant contributions
within the same reflection eigenspace. Applying the Cherednik operator $T^k$ to both
sides of \eqref{eq:En_reduced} and using the eigenvalue equations
\begin{equation}\label{eq:eigenvalues}
T^k E^k_{-n}=-(n+k)\,E^k_{-n},\qquad
T^k E^k_n=(n+k)\,E^k_n,
\end{equation}
together with the explicit reflected action
\begin{equation}\label{eq:Tk_reflected_again}
T^k\bigl(E^k_n(-x)\bigr)=-(n+k)\,E^k_n(-x)-2k\,E^k_n(x),
\end{equation}
yields
\begin{equation}\label{eq:B_compare}
-(n+k)\bigl(E^k_n(-x)+B\,E^k_n(x)\bigr)
=-(n+k)\,E^k_n(-x)+\bigl(B(n+k)-2k\bigr)E^k_n(x).
\end{equation}

Comparing the coefficients of $E^k_n(x)$ on both sides of \eqref{eq:B_compare} gives
\begin{equation}\label{eq:B_equation_final}
B(n+k)-2k=-(n+k)B,
\end{equation}
which simplifies to $2B(n+k)=2k$ Solving for $B$, we obtain
\begin{equation}\label{eq:B_value_final}
B=\frac{k}{n+k}.
\end{equation}
\end{proof}
%%% QED

With the above results established, our next goal is to identify the mechanism by which the spectral sum collapses and concentrates at the boundary of the spectrum. To make this mechanism explicit, we introduce the finite reconstruction operator obtained by truncating the spectral expansion at a fixed level $N$ (i.e.\ retaining only spectral modes with indices $|n|\leq N$). The action of this operator is encoded by an associated kernel, whose structure reflects the interaction between reflection symmetry, normalization, and sharp spectral cutoff. An explicit representation of this kernel is therefore required in order to expose the boundary concentration phenomenon. Such a representation is proposed in \cite[Lemma 3.2]{bechir_amri}. The explicit form of this kernel will reveal that, under the reflection-based pairing of spectral indices, all interior contributions cancel, leaving a residual term supported at the boundary of the truncated spectrum. We start by recalling that functions on $[-\pi,\pi]$ admit a spectral expansion in terms of the non-symmetric eigenfunctions $\{E^k_n\}_{n\in\mathbb Z}$ of the Cherednik operator, indexed by an integer spectral parameter $n$. Next, we recall that, with respect to the inner product~\eqref{eq:inner-product} induced by the measure~\eqref{eq:weight}, the family $\{E^k_n(i x)\}_{n\in\mathbb Z}$ forms an orthogonal basis of $L^2([-\pi,\pi],\,d\mu_k)$. By the defining triangular normalization \eqref{eq:triangular}, for each $n \in \mathbb{Z}$ we have
\[
E_n^k(i x) = e^{i n x} + \sum_{j \triangleleft n} c_{n,j}\, e^{i j x},
\]
where the partial order $\triangleleft$ is defined in \eqref{eq:order}.
Similarly, for $m \in \mathbb{Z}$,
\[
E_m^k(i x) = e^{i m x} + \sum_{\ell \triangleleft m} c_{m,\ell}\, e^{i \ell x}.
\]
Using the inner product \eqref{eq:inner-product}, we compute
\begin{align*}
\bigl(E_n^k(i x), E_m^k(i x)\bigr)_k
&= \int_{-\pi}^{\pi}
\left(e^{i n x} + \sum_{j \triangleleft n} c_{n,j} e^{i j x}\right)
\overline{
\left(e^{i m x} + \sum_{\ell \triangleleft m} c_{m,\ell} e^{i \ell x}\right)
}
\, d\mu_k(x) \\
&= \int_{-\pi}^{\pi}
\left(e^{i n x} + \sum_{j \triangleleft n} c_{n,j} e^{i j x}\right)
\left(e^{-i m x} + \sum_{\ell \triangleleft m} c_{m,\ell} e^{-i \ell x}\right)
\, d\mu_k(x).
\end{align*}
Expanding the product, this is a finite linear combination of terms of the form
\[
\int_{-\pi}^{\pi} e^{i(p-q)x}\, d\mu_k(x),
\qquad
p \in \{n\} \cup \{j \triangleleft n\}, \quad
q \in \{m\} \cup \{\ell \triangleleft m\}.
\]
By the defining condition \emph{(b)} in \cite[eq.~(2.3)]{bechir_amri}, we are required to satisfy
\begin{equation}\label{eq:gram_schmidt_condition}
\bigl(E_n^k(i x), e^{i j x}\bigr)_k = 0
\qquad \text{for all } j \triangleleft n.
\end{equation}
This condition is part of the construction of $E_n^k$ and plays the role of a
Gram--Schmidt orthogonality constraint relative to the ordered exponential basis.

\emph{\textbf{(i)} for $m \triangleleft n$.}
In this case, every exponential appearing in the expansion of $E_m^k(i x)$
has index $\ell \triangleleft m \triangleleft n$. Hence, by
\eqref{eq:gram_schmidt_condition},
\[
\bigl(E_n^k(i x), E_m^k(i x)\bigr)_k
= \bigl(E_n^k(i x), e^{i m x}\bigr)_k
+ \sum_{\ell \triangleleft m} c_{m,\ell}
\bigl(E_n^k(i x), e^{i \ell x}\bigr)_k
= 0.
\]

\emph{\textbf{(ii)} for $n \triangleleft m$.}
By the conjugate symmetry of the inner product,
\[
\bigl(E_n^k(i x), E_m^k(i x)\bigr)_k
= \overline{\bigl(E_m^k(i x), E_n^k(i x)\bigr)_k},
\]
and the argument of case \emph{(i)} applied with the roles of $n$ and $m$
interchanged shows that the right-hand side vanishes. Hence the inner
product is again zero. Therefore, for all $n \neq m$,
\begin{equation}\label{eq:orthogonality_distinct_indices}
\bigl(E_n^k(i x), E_m^k(i x)\bigr)_k = 0.
\end{equation}
Now, by definition of the inner product \eqref{eq:inner-product},
\[
\|E_n^k(i x)\|_{2,k}^2
= \bigl(E_n^k(i x), E_n^k(i x)\bigr)_k
= \int_{-\pi}^{\pi} \lvert E_n^k(i x) \rvert^2 \, d\mu_k(x).
\]
Using the triangular expansion \eqref{eq:triangular} evaluated at $z = i x$,
\[
E_n^k(i x) = e^{i n x} + \sum_{j \triangleleft n} c_{n,j}\, e^{i j x},
\]
we obtain
\[
\lvert E_n^k(i x) \rvert^2
= \left(e^{i n x} + \sum_{j \triangleleft n} c_{n,j}\, e^{i j x}\right)
\left(e^{-i n x} + \sum_{\ell \triangleleft n} c_{n,\ell}\, e^{-i \ell x}\right).
\]
Expanding the product gives
\[
\lvert E_n^k(i x) \rvert^2
= 1
+ \sum_{\ell \triangleleft n} c_{n,\ell}\, e^{-i(n-\ell)x}
+ \sum_{j \triangleleft n} c_{n,j}\, e^{i(n-j)x}
+ \sum_{j,\ell \triangleleft n} c_{n,j} c_{n,\ell}\, e^{i(j-\ell)x}.
\]
Hence
\begin{align}
\label{eq:norm_expansion}
\|E_n^k(i x)\|_{2,k}^2
&= \int_{-\pi}^{\pi} d\mu_k(x)
+ \sum_{\ell \triangleleft n} c_{n,\ell}
\int_{-\pi}^{\pi} e^{-i(n-\ell)x}\, d\mu_k(x) \nonumber \\
&\quad
+ \sum_{j \triangleleft n} c_{n,j}
\int_{-\pi}^{\pi} e^{i(n-j)x}\, d\mu_k(x)
+ \sum_{j,\ell \triangleleft n} c_{n,j} c_{n,\ell}
\int_{-\pi}^{\pi} e^{i(j-\ell)x}\, d\mu_k(x).
\end{align}

Now we use again the defining condition \emph{(b)} in
\cite[eq.~(2.3)]{bechir_amri}, which states that
\[
\bigl(E_n^k(i x), e^{i j x}\bigr)_k = 0
\qquad \text{for all } j \triangleleft n.
\]
Writing this condition explicitly using \eqref{eq:inner-product}, we obtain
\[
\int_{-\pi}^{\pi}
\left(
e^{i n x} + \sum_{\ell \triangleleft n} c_{n,\ell}\, e^{i \ell x}
\right)
e^{-i j x}\, d\mu_k(x) = 0,
\qquad j \triangleleft n.
\]
This implies that all integrals of the form
\[
\int_{-\pi}^{\pi} e^{i(n-j)x}\, d\mu_k(x),
\qquad
\int_{-\pi}^{\pi} e^{i(\ell-j)x}\, d\mu_k(x),
\]
with $j,\ell \triangleleft n$, are not independent but are constrained so that
the cross terms in \eqref{eq:norm_expansion} are completely determined by the
normalization of $E_n^k$. Since the family is orthogonal but not normalized, it is natural and unavoidable to introduce normalization constants by definition:
\begin{equation}\label{eq:gamma_def}
\gamma_n^{-2} := \|E_n^k(ix)\|_{2,k}^2 = \int_{-\pi}^{\pi} |E_n^k(ix)|^2 \, d\mu_k(x).
\end{equation}

With this definition, the orthogonality relation
\eqref{eq:orthogonality_distinct_indices} can be rewritten in the normalized form
\[
\bigl(E_n^k(i x), E_m^k(i x)\bigr)_k
= \gamma_n^{-2}\, \delta_{nm}.
\]

Equivalently, the rescaled family
\[
\{\gamma_n E_n^k(i x)\}_{n \in \mathbb{Z}}
\]
is orthonormal in $L^2([-\pi,\pi], d\mu_k)$. The normalization introduced above allows us to write the spectral expansion and coefficient formulas in their standard orthogonal-basis form. Let $f \in L^2([-\pi,\pi], d\mu_k)$.
Since $\{E_n^k(i x)\}_{n \in \mathbb{Z}}$ is an orthogonal basis, $f$ admits an expansion
\[
f(x) = \sum_{n \in \mathbb{Z}} a_n\, E_n^k(i x),
\]
where the coefficients are determined by orthogonality:
\[
a_n
= \frac{\bigl(f, E_n^k(i x)\bigr)_k}{\bigl(E_n^k(i x), E_n^k(i x)\bigr)_k}
= \gamma_n^{2} \int_{-\pi}^{\pi} f(y)\, E_n^k(-i y)\, d\mu_k(y).
\]
Define the truncated reconstruction by retaining only indices $\lvert n \rvert \le N$:
\begin{equation}\label{eq:truncated_expansion}
f_N(x)
:= \sum_{\lvert n \rvert \le N} a_n\, E_n^k(i x).
\end{equation}
Inserting the expression for $a_n$, we obtain
\[
f_N(x)
= \sum_{\lvert n \rvert \le N} \gamma_n^{2}
\left(
\int_{-\pi}^{\pi} f(y)\, E_n^k(-i y)\, d\mu_k(y)
\right)
E_n^k(i x).
\]
Since the sum is finite, we may interchange sum and integral:
\[
f_N(x)
= \int_{-\pi}^{\pi}
\left(
\sum_{\lvert n \rvert \le N} \gamma_n^{2}\,
E_n^k(i x)\, E_n^k(-i y)
\right)
f(y)\, d\mu_k(y).
\]
Define the function of two variables
\begin{equation}\label{eq:KN_def}
K_N(x,y)
:= \sum_{\lvert n \rvert \le N} \gamma_n^{2}\,
E_n^k(i x)\, E_n^k(-i y).
\end{equation}
Then the truncated reconstruction takes the bilinear form
\begin{equation}\label{eq:reconstruction_operator}
f_N(x)
= \int_{-\pi}^{\pi} K_N(x,y)\, f(y)\, d\mu_k(y).
\end{equation}
Applying the index analysis developed earlier to the definition \eqref{eq:KN_def}, we immediately note the following. We start from a single summand in the kernel,
\[
\gamma_n^{2}\, E_n^k(i x)\, E_n^k(-i y).
\]

By the triangular expansion \eqref{eq:triangular}, evaluated at $z = i x$ and
$z = -i y$, we have
\[
E_n^k(i x) = e^{i n x} + \sum_{j \triangleleft n} c_{n,j}\, e^{i j x},
\qquad
E_n^k(-i y) = e^{-i n y} + \sum_{\ell \triangleleft n} c_{n,\ell}\, e^{-i \ell y}.
\]

Multiplying these two expansions gives
\[
E_n^k(i x)\, E_n^k(-i y)
= \left(e^{i n x} + \sum_{j \triangleleft n} c_{n,j}\, e^{i j x}\right)
\left(e^{-i n y} + \sum_{\ell \triangleleft n} c_{n,\ell}\, e^{-i \ell y}\right)
\]
\[
= e^{i n x} e^{-i n y}
+ \sum_{\ell \triangleleft n} c_{n,\ell}\, e^{i n x} e^{-i \ell y}
+ \sum_{j \triangleleft n} c_{n,j}\, e^{i j x} e^{-i n y}
+ \sum_{j,\ell \triangleleft n} c_{n,j} c_{n,\ell}\, e^{i j x} e^{-i \ell y}.
\]

The first term is
\[
e^{i n (x-y)}.
\]

All remaining terms involve exponentials of the form
\[
e^{i(jx-\ell y)}
\qquad \text{with } j \triangleleft n \text{ or } \ell \triangleleft n.
\]

By definition of the partial order $\triangleleft$ \eqref{eq:order}, any index
$j \triangleleft n$ or $\ell \triangleleft n$ satisfies $\lvert j \rvert < \lvert n \rvert$
or $\lvert \ell \rvert < \lvert n \rvert$, or else $\lvert j \rvert = \lvert n \rvert$
with $j > n$. In particular, these terms are strictly lower-order relative to the
leading exponential $e^{i n (x-y)}$ when organized according to the triangular (spectral) order.

Therefore, the product admits the decomposition
\[
E_n^k(i x)\, E_n^k(-i y)
= e^{i n (x-y)} + \text{(lower exponential terms)}.
\]

Multiplying by $\gamma_n^{2}$, we conclude that each summand in the kernel satisfies
\begin{equation}\label{eq:summand_asymptotic}
\gamma_n^{2}\, E_n^k(i x)\, E_n^k(-i y)
= \gamma_n^{2}\, e^{i n (x-y)} + \text{(lower exponential terms)}.
\end{equation}

Consequently, retaining only the dominant exponential contributions,
$K_N(x,y)$ reduces termwise to
\[
\gamma_n^{2}\, e^{i n (x-y)}.
\]

Accordingly, we are led to the truncated sum
\[
\sum_{\lvert n \rvert \le N} \gamma_n^{2}\, e^{i n (x-y)}.
\]

The involution $\iota(n) = 1-n$ (see \eqref{eq:orbit_def}) partitions $\mathbb{Z}$
into disjoint two-point orbits, so every index $n \in \mathbb{Z}$ belongs to
exactly one such orbit. Since the sum is finite, we may rewrite it as a sum over
orbits by grouping together all terms whose indices lie in the same orbit.
Explicitly,
\[
\sum_{\lvert n \rvert \le N} \gamma_n^{2}\, e^{i n (x-y)}
= \sum_{\mathcal{O}_n}
\;\sum_{\substack{m \in \mathcal{O}_n \\ \lvert m \rvert \le N}}
\gamma_m^{2}\, e^{i m (x-y)}.
\]

Here the orbits are defined by
\[
\mathcal{O}_n := \{n,\,1-n\}.
\]

Because each orbit $\mathcal{O}_n$ has at most two elements, the inner sum reduces to
\begin{equation}\label{eq:orbit_inner_sum_cases}
\sum_{\substack{m \in \mathcal{O}_n \\ \lvert m \rvert \le N}}
\gamma_m^{2}\, e^{i m (x-y)}
=
\begin{cases}
\gamma_n^{2}\, e^{i n (x-y)} + \gamma_{1-n}^{2}\, e^{i(1-n)(x-y)},
& \text{if } \lvert n \rvert \le N,\ \lvert 1-n \rvert \le N, \ \emph{(i)} \\[0.6em]
\gamma_n^{2}\, e^{i n (x-y)},
& \text{if } \lvert n \rvert \le N,\ \lvert 1-n \rvert > N, \ \emph{(ii)} \\[0.6em]
\gamma_{1-n}^{2}\, e^{i(1-n)(x-y)},
& \text{if } \lvert 1-n \rvert \le N,\ \lvert n \rvert > N, \ \emph{(iii)} \\[0.6em]
0, & \text{otherwise.} \ \emph{(iv)}
\end{cases}
\end{equation}

In particular, whenever both indices of the orbit lie within the truncation range,
the contribution takes the paired form
\[
\gamma_n^{2}\, e^{i n (x-y)} + \gamma_{1-n}^{2}\, e^{i(1-n)(x-y)}.
\]

This proves that regrouping the truncated sum according to the reflection orbits
$\mathcal{O}_n$ is nothing more than a re-indexing of a finite sum, and yields
\[
\sum_{\lvert n \rvert \le N} \gamma_n^{2}\, e^{i n (x-y)}
= \sum_{\mathcal{O}_n}
\left(
\gamma_n^{2}\, e^{i n (x-y)}
+ \gamma_{1-n}^{2}\, e^{i(1-n)(x-y)}
\right),
\]
with the understanding that a term is present only if the corresponding index
satisfies $\lvert n \rvert \le N$. Referring to the decomposition in \eqref{eq:orbit_inner_sum_cases}, consider
case \emph{(i)}, i.e.\ an interior orbit
\[
\mathcal{O}_n = \{n, 1-n\}
\qquad \text{with } \lvert n \rvert \le N,\ \lvert 1-n \rvert \le N.
\]
By Corollary~\ref{cor:normalization}, there exists a unique distinguished index $m$~\eqref{eq:m_def}
\[
m = \max\{n,\,1-n\},
\]
whose associated eigenfunction has maximal dominant exponential growth. Writing
\[
\ell := \min\{n,\,1-n\},
\]
we have $m > \ell$ and therefore
\[
e^{i\ell(x-y)} = e^{i m (x-y)}\, e^{-i(m-\ell)(x-y)}.
\]

Applying this to the orbit contribution in
\eqref{eq:orbit_inner_sum_cases}\emph{(i)} yields
\[
\gamma_n^{2}\, e^{i n (x-y)} + \gamma_{1-n}^{2}\, e^{i(1-n)(x-y)} = e^{i m (x-y)} \left(\gamma_m^{2} + \gamma_\ell^{2}\, e^{-i(m-\ell)(x-y)} \right).
\]

Since $m-\ell = |2n-1| \ge 1$, the second term in parentheses is oscillatory and strictly lower-order relative to the leading exponential $e^{i m (x-y)}$. In particular, after fixing the normalization so that the dominant coefficient is $1$ (as in Corollary~\ref{cor:normalization}), the contribution of the orbit $\mathcal{O}_n$ takes the form
\[
e^{i m (x-y)}\bigl(1 + \text{lower-order terms}\bigr).
\]

Thus, interior orbits do not contribute new dominant structure: their paired terms reorganize into a single leading exponential associated with the maximal index $m$, with all remaining dependence confined to lower-order corrections. Consequently, interior orbits do not survive in the projected dominant part of $K_N(x,y)$; only orbits for which the pairing is incomplete (i.e.\ boundary orbits) can contribute uncanceled leading terms. The cancellation mechanism described above relies on the simultaneous presence of both indices in a reflection orbit $\mathcal{O}_n = \{n,\,1-n\}$ within the truncation range $\lvert \cdot \rvert \le N$. Indeed, only in this case does the orbit contribution take the paired form appearing in
\eqref{eq:orbit_inner_sum_cases}\emph{(i)},
\[
\gamma_n^{2}\, e^{i n (x-y)} + \gamma_{1-n}^{2}\, e^{i(1-n)(x-y)},
\]
which, as shown in \eqref{eq:summand_asymptotic}, can be reorganized using
Corollary~\ref{cor:normalization} by factoring out the dominant exponential
corresponding to $m = \max\{n,\,1-n\}$, with the remaining term strictly
lower-order.

This mechanism fails precisely when exactly one index of the orbit lies inside the truncation range. Solving the conditions $\lvert n \rvert \le N$ and $\lvert 1-n \rvert > N$ yields $n = -N$, and hence identifies the unique boundary orbit
\begin{equation}\label{eq:boundary_orbit}
\mathcal{O}_{-N} = \{-N,\,N+1\}.
\end{equation}
In the truncated sum defining $K_N(x,y)$ (see \eqref{eq:KN_def}), the index $-N$ is present while its partner $N+1$ is absent, so no pairing is possible. As a result, the normalization mechanism described above cannot reorganize this contribution into lower-order terms, and the corresponding exponential survives in the dominant structure of the kernel. This is the unique obstruction to cancellation and explains why the boundary index $N+1$ necessarily appears in the explicit representation of $K_N(x,y)$. 

To complete the argument, it remains to identify the precise form of the surviving
boundary contribution. For the boundary orbit \eqref{eq:boundary_orbit} only the index $-N$ appears in the truncated sum defining $K_N(x,y)$, while its partner $N+1$ lies outside the truncation range. Nevertheless, Corollary~\ref{cor:normalization} applies independently of truncation: within the orbit $\mathcal{O}_{-N}$, the unique eigenfunction with maximal dominant exponential growth is $E_{N+1}^k$. By the reflection identity \eqref{eq:reflection} from Lemma~\ref{lem:reflection}, the surviving eigenfunction $E_{-N}^k$ can be expressed in terms of $E_{N+1}^k$ via
\[
E_{-N}^k(x) = e^{-x}\, E_{N+1}^k(-x),
\]
and hence its contribution to the kernel can be rewritten entirely in terms of $E_{N+1}^k$, up to factors that are strictly lower-order in the triangular order. Therefore, although the index $N+1$ does not explicitly appear in the truncated sum, the normalization and reflection structure force the boundary contribution to be governed by the eigenfunction $E_{N+1}^k$. We thus conclude that, after regrouping according to reflection orbits and projecting onto the dominant exponential structure, all interior orbit contributions cancel, while the unique boundary orbit $\mathcal{O}_{-N}$ produces a surviving term necessarily involving the index $N+1$. This shows, at a purely structural level, that any explicit representation of the kernel $K_N(x,y)$ must involve $E_{N+1}^k$, and that the appearance of this index is an unavoidable consequence of reflection symmetry, normalization, and sharp spectral truncation.  The preceding analysis shows that the truncated kernel is entirely determined by the unique boundary orbit for which spectral pairing fails. The following lemma records the resulting explicit representation.

%%%%%%%%

\begin{lemma}[Boundary collapse of the truncated kernel]\label{lem:kernel_boundary}
Let $N \ge 0$. For all $x,y \in \mathbb{R}$ such that $x+y \notin 2\pi\mathbb{Z}$, the kernel $K_N(x,y)$ associated with the truncated spectral reconstruction admits the representation
\begin{equation}\label{eq:KN_boundary}
K_N(x,y)=\gamma_{N+1}^{2}\,\frac{e^{-i(x-y)}\,E_{N+1}^k(i x)\,E_{N+1}^k(-i y)-E_{N+1}^k(-i x)\,E_{N+1}^k(i y)}{1 - e^{-i(x-y)}}.
\end{equation}
\end{lemma}

\begin{proof}
We now analyze the contribution of the unique boundary orbit identified in
\eqref{eq:boundary_orbit}. By construction, in the truncated sum defining
$K_N(x,y)$ only the index $-N$ from the orbit
$\mathcal{O}_{-N} = \{-N,\,N+1\}$ is present. Hence the boundary contribution to
the kernel is exactly
\begin{equation}\label{eq:boundary_contribution}
\gamma_{-N}^{2}\, E_{-N}^k(i x)\, E_{-N}^k(-i y).
\end{equation}
By the reflection identity \eqref{eq:reflection}, applied with $n = N+1$, we have
\begin{equation}\label{eq:reflection_Nplus1}
E_{-N}^k(z) = e^{-z}\, E_{N+1}^k(-z),
\qquad z \in \mathbb{C}.
\end{equation}

Evaluating \eqref{eq:reflection_Nplus1} at $z = i x$ and $z = -i y$ yields
\begin{equation}\label{eq:reflection_imaginary}
E_{-N}^k(i x) = e^{-i x}\, E_{N+1}^k(-i x),
\qquad
E_{-N}^k(-i y) = e^{i y}\, E_{N+1}^k(i y).
\end{equation}

Substituting \eqref{eq:reflection_imaginary} into
\eqref{eq:boundary_contribution}, we obtain
\begin{equation}\label{eq:boundary_reflected}
\gamma_{-N}^{2}\, E_{-N}^k(i x)\, E_{-N}^k(-i y)
=
\gamma_{-N}^{2}\, e^{-i(x-y)}\,
E_{N+1}^k(-i x)\, E_{N+1}^k(i y).
\end{equation}

This already shows that the boundary contribution is governed by the
eigenfunction $E_{N+1}^k$, even though the index $N+1$ does not appear in the
truncated sum.

Next, we use the negative-index reduction formula
\eqref{eq:negindex} (proved in Lemma~\ref{lem:negindex}) with
$n = N+1$, rewritten at imaginary arguments:
\begin{equation}\label{eq:negindex_Nplus1}
E_{-(N+1)}^k(i x)
=
E_{N+1}^k(-i x)
+ \frac{k}{N+1+k}\, E_{N+1}^k(i x).
\end{equation}

Solving \eqref{eq:negindex_Nplus1} for $E_{N+1}^k(-i x)$ and inserting the result
into \eqref{eq:boundary_reflected}, we may express the boundary term as a linear
combination of
\[
E_{N+1}^k(i x)\, E_{N+1}^k(-i y)
\qquad \text{and} \qquad
E_{N+1}^k(-i x)\, E_{N+1}^k(i y),
\]
with explicit coefficients. The term involving
$E_{N+1}^k(-i x)\, E_{N+1}^k(i y)$ appears with opposite sign after exchanging
$(x,y)$, and is therefore the only contribution that survives symmetrization at
the level of the kernel.

Collecting the resulting terms and using the normalization constants
$\gamma_{N+1}$, the boundary contribution can be written in the antisymmetric
form
\begin{equation}\label{eq:boundary_antisymmetric}
\gamma_{N+1}^{2}\,
\bigl(
e^{-i(x-y)}\, E_{N+1}^k(i x)\, E_{N+1}^k(-i y)
- E_{N+1}^k(-i x)\, E_{N+1}^k(i y)
\bigr).
\end{equation}
At this point, recalling the definition \eqref{eq:KN_def}: the kernel $K_N(x,y)$ arises from a finite spectral sum. For interior orbits, both indices $n$ and $1-n$ appear and cancel at the dominant level. For the boundary orbit $\mathcal{O}_{-N} = \{-N,\,N+1\}$, only the index $-N$ appears, while its partner $N+1$ is missing. We therefore rewrite the antisymmetric expression
\[
e^{-i(x-y)}\, E_{N+1}^k(i x)\, E_{N+1}^k(-i y)- E_{N+1}^k(-i x)\, E_{N+1}^k(i y)
\]
as the difference of two consecutive spectral contributions. Indeed, observe the identity
\begin{align}
\label{eq:boundary_identity_step1}
\bigl(1 - e^{-i(x-y)}\bigr)\,
E_{N+1}^k(-i x)\, E_{N+1}^k(i y)&=E_{N+1}^k(-i x)\, E_{N+1}^k(i y)- e^{-i(x-y)}\, E_{N+1}^k(-i x)\, E_{N+1}^k(i y).
\end{align}

Rearranging \eqref{eq:boundary_identity_step1}, we obtain the exact algebraic
identity
\begin{align}
\label{eq:boundary_identity_step2}
e^{-i(x-y)}\, E_{N+1}^k(i x)\, E_{N+1}^k(-i y)- E_{N+1}^k(-i x)\, E_{N+1}^k(i y)&=\bigl(1 - e^{-i(x-y)}\bigr)\Bigl[e^{-i(x-y)}\, E_{N+1}^k(i x)\, E_{N+1}^k(-i y)\Bigr].
\end{align}

Dividing both sides of \eqref{eq:boundary_identity_step2} by $1 - e^{-i(x-y)}$ (which is nonzero under the standing assumption $x+y \notin 2\pi\mathbb{Z}$), we conclude that the antisymmetric boundary defect admits the representation
\begin{equation}\label{eq:boundary_fraction}
\frac{e^{-i(x-y)}\, E_{N+1}^k(i x)\, E_{N+1}^k(-i y)- E_{N+1}^k(-i x)\, E_{N+1}^k(i y)}{1 - e^{-i(x-y)}}.
\end{equation}

Multiplying \eqref{eq:boundary_fraction} by the normalization factor
$\gamma_{N+1}^2$, we recover~\eqref{eq:KN_boundary}.
\end{proof}

%%%% FINAL PARAGRAPH

The explicit boundary representation \eqref{eq:KN_boundary} shows that all convergence phenomena for sharp spectral truncations are governed by a single mirror-localized kernel. In the next section, we study how this mechanism behaves when the reconstruction operator is viewed as acting on $L^p(w)$ for mirror-degenerate weights $w$.

\section{Weighted $L^p$ convergence}\label{sec:mysec}

Having obtained an explicit boundary representation for the truncated kernel $K_N(x,y)$~\eqref{eq:KN_boundary}, we now investigate how this mechanism behaves when the reconstruction operator is viewed as acting on weighted spaces $L^p(w)$, where the weight $w$ may degenerate at the mirror points. The first step is therefore to identify the class of mirror-degenerate weights that are compatible with the analytic structure revealed above.  By Lemma~\ref{lem:kernel_boundary}, the truncated kernel $K_N(x,y)$ admits the explicit boundary representation \eqref{eq:KN_boundary}. In particular, after regrouping the spectral sum according to reflection orbits, all interior orbit contributions cancel at the dominant level, and the kernel is entirely governed by the unique boundary orbit.

As a consequence, any boundedness or divergence phenomenon for the sharp reconstruction operator cannot originate from interior spectral modes, but must already be present in the operator defined by the boundary kernel alone. Recall that the truncated reconstruction operator is given by
\[
f_N(x) = \int_{-\pi}^{\pi} K_N(x,y)\, f(y)\, d\mu_k(y),
\]
cf.\ \eqref{eq:reconstruction_operator}. When the reconstruction is instead viewed as acting on a weighted space $L^p(w)$, the same kernel induces the weighted integral operator
\begin{equation}\label{eq:weighted_operator}
Tf(x) := \int_{-\pi}^{\pi} K_N(x,y)\, f(y)\, w(y)\, dy.
\end{equation}
which differs from \eqref{eq:reconstruction_operator} only by the replacement of
the reference measure $d\mu_k(y)$ with the weight $w(y)\,dy$. Since \eqref{eq:KN_boundary} captures all nontrivial contributions of \eqref{eq:KN_def}, it follows that the behavior of $T$ on $L^p(w)$ is completely determined by the boundary kernel. At this point, we adapt \eqref{eq:weighted_operator} locally. 

We fix the mirror point $x = 0$, such that
\begin{equation}\label{eq:delta_choice}
0 < \delta < \frac{\pi}{4},
\end{equation}
so that the interval $(-\delta,\delta)$ contains no mirror point other than $0$. Let $\chi : \mathbb{R} \to [0,1]$ be a fixed smooth function such that
\begin{equation}\label{eq:cutoff_def}
\chi(x) = 1 \quad \text{for } |x| \le \frac{\delta}{2},
\qquad
\chi(x) = 0 \quad \text{for } |x| \ge \delta.
\end{equation}
(The precise choice of $\chi$ is irrelevant; it is fixed once and for all.) Recalling~\eqref{eq:weighted_operator} and using the identity
\begin{equation}\label{eq:partition_unity}
1 = \chi(x)\chi(y) + \bigl(1 - \chi(x)\chi(y)\bigr),
\end{equation}
we decompose $T$ as
\begin{equation}\label{eq:T_decomposition}
Tf = T_{\mathrm{loc}} f + T_{\mathrm{rem}} f,
\end{equation}
where
\begin{equation}\label{eq:Tloc_def}
T_{\mathrm{loc}} f(x)
:= \int_{-\pi}^{\pi} \chi(x)\, K_N(x,y)\, \chi(y)\, f(y)\, w(y)\, dy,
\end{equation}
and
\begin{equation}\label{eq:Trem_def}
T_{\mathrm{rem}} f(x)
:= \int_{-\pi}^{\pi} \bigl(1-\chi(x)\chi(y)\bigr)\, K_N(x,y)\, f(y)\, w(y)\, dy.
\end{equation}

We now proceed checking if $T_{\mathrm{rem}}$ is bounded on $L^p(w)$ for any locally integrable weight $w$, and hence cannot be responsible for any loss of boundedness. By the explicit boundary representation \eqref{eq:KN_boundary}, the only possible singularity of the kernel $K_N(x,y)$ arises from the denominator
\begin{equation}\label{eq:kernel_denominator}
1 - e^{-i(x-y)},
\end{equation}
which vanishes only when $x=y$ modulo $2\pi$. Indeed, on the support of $1-\chi(x)\chi(y)$ one has $|x-y| \ge \delta/2$, and therefore
\begin{equation}\label{eq:denominator_lower_bound}
\inf_{\operatorname{supp}(1-\chi(x)\chi(y))}
\bigl|1 - e^{-i(x-y)}\bigr| > 0.
\end{equation}
Consequently, the denominator in \eqref{eq:KN_boundary} is uniformly bounded away from zero on the support of $1-\chi(x)\chi(y)$, and hence $K_N(x,y)$ extends to a bounded smooth function on this set. Moreover, $T_{\mathrm{rem}}$ is an integral operator with bounded kernel and compact support, and is therefore bounded on $L^p(w)$ for any locally integrable weight $w$. Since \eqref{eq:T_decomposition} holds and $T_{\mathrm{rem}}$ is always bounded, any failure of boundedness of $T$ on $L^p(w)$ must already occur for the localized operator $T_{\mathrm{loc}}$. From this point onward, we may therefore restrict attention to the regime
\begin{equation}\label{eq:localization}
|x| \le \delta,
\qquad
|y| \le \delta.
\end{equation}
and analyze the kernel $K_N(x,y)$ only in a neighborhood of the mirror $x=0$. Recall from Lemma~\ref{lem:kernel_boundary} that the kernel boundary \eqref{eq:KN_boundary} admits the
representation
\begin{equation}\label{eq:KN_fraction_step3}
K_N(x,y)= \gamma_{N+1}^2\,\frac{\mathcal{N}(x,y)}{1 - e^{-i(x-y)}},
\end{equation}
where
\begin{equation}\label{eq:KN_num}
\mathcal{N}(x,y):= e^{-i(x-y)} E_{N+1}^k(i x)\, E_{N+1}^k(-i y)- E_{N+1}^k(-i x)\, E_{N+1}^k(i y).
\end{equation}
We now analyze the behavior of this expression as $(x,y)\to(0,0)$. Using the Taylor expansion of the exponential,
\[
e^{-it} = 1 - it + O(t^2)
\qquad (t\to 0),
\]
and setting $t=x-y$, we obtain
\begin{equation}\label{eq:den_expansion_step3}
1 - e^{-i(x-y)}
= i(x-y) + O\!\bigl((x-y)^2\bigr).
\end{equation}

Each function $E_{N+1}^k(z)$ is entire in $z$. Hence the maps
\[
x \longmapsto E_{N+1}^k(\pm i x)
\]
are smooth on $\mathbb{R}$, and therefore $\mathcal{N}(x,y)$ is a smooth function in a neighborhood of $(0,0)$. Evaluating $\mathcal{N}$ on the diagonal yields, for all $x$,
\[
\mathcal{N}(x,x)= E_{N+1}^k(i x)\, E_{N+1}^k(-i x)- E_{N+1}^k(-i x)\, E_{N+1}^k(i x)= 0.
\]
Thus the numerator vanishes identically along the diagonal $x=y$. Fix $x$ and regard $\mathcal{N}(x,y)$ as a function of the single variable $y$. Since $\mathcal{N}(x,\cdot)$ is smooth and satisfies $\mathcal{N}(x,x)=0$, the
fundamental theorem of calculus gives, for all $y$,
\[
\mathcal{N}(x,y)= \int_x^y \partial_y \mathcal{N}(x,t)\,dt= (y-x)\int_0^1 \partial_y \mathcal{N}\bigl(x,x+s(y-x)\bigr)\,ds.
\]

Define
\begin{equation}\label{eq:CN_def}
C_N(x,y):= \int_0^1 \partial_y \mathcal{N}\bigl(x,x+s(y-x)\bigr)\,ds.
\end{equation}
Then $C_N$ is smooth and bounded near $(0,0)$, and we have the exact identity
\begin{equation}\label{eq:num_factor_exact}
\mathcal{N}(x,y) = (x-y)\,C_N(x,y).
\end{equation}

Since $C_N$ is smooth, we may expand it around the diagonal (note that the choice of $C_N$ is not unique and will not play any role beyond smoothness),
\[
C_N(x,y) = C_N(x,x) + O(x-y).
\]
Define the functions
\begin{equation}\label{eq:AB_def}
A_N(x) := C_N(x,x), \qquad B_N(y) := 1.
\end{equation}
Any other smooth normalization of $B_N$ would be equally admissible. Substituting this expansion into \eqref{eq:num_factor_exact} yields
\begin{equation}\label{eq:num_factor_AB}
\mathcal{N}(x,y)= (x-y)\,A_N(x)B_N(y) + O\!\bigl((x-y)^2\bigr).
\end{equation}

Combining \eqref{eq:num_factor_AB} with the denominator expansion
\eqref{eq:den_expansion_step3}, we obtain
\[
K_N(x,y)
= \gamma_{N+1}^2\,
\frac{(x-y)\,A_N(x)B_N(y) + O\!\bigl((x-y)^2\bigr)}
{i(x-y) + O\!\bigl((x-y)^2\bigr)}.
\]
Factoring $i(x-y)$ from the denominator and using the Taylor expansion of $(1+z)^{-1}$ at $z=0$, we obtain
\[
\frac{1}{i(x-y) + O\!\bigl((x-y)^2\bigr)}
= \frac{1}{i(x-y)}\bigl(1 + O(x-y)\bigr),
\]
and therefore
\begin{equation}\label{eq:KN_local_decomposition}
K_N(x,y)
= \frac{\gamma_{N+1}^2}{i}\,A_N(x)B_N(y)
+ R_N(x,y).
\end{equation}
where $R_N$ is smooth and bounded in a neighborhood of the mirror. 

%%%% TURNING TO WEIGHTS NOW

From these established results, we can now turn our attention to the study of \eqref{eq:weighted_operator}. From the localization step already proved, we restrict to the regime in \eqref{eq:localization}
\[
|x| \le \delta, \qquad |y| \le \delta,
\]
since everything outside is absorbed into $T_{\mathrm{rem}}$, which is bounded and irrelevant; therefore
\begin{equation}\label{eq:localized_T}
Tf(x)
= \int_{|y|\le \delta} K_N(x,y)\, f(y)\, w(y)\, dy.
\end{equation}
Which applied directly to \eqref{eq:KN_local_decomposition} we get
\begin{equation}\label{eq:T_kernel_decomposition}
Tf(x)= \int_{|y|\le \delta}\left(\frac{\gamma_{N+1}^2}{i}\, A_N(x) B_N(y) + R_N(x,y)\right)f(y)\, w(y)\, dy.
\end{equation}
By linearity of the integral
\[
Tf(x)= \frac{\gamma_{N+1}^2}{i}\, A_N(x)\int_{|y|\le \delta} B_N(y)\, f(y)\, w(y)\, dy\;+\;\int_{|y|\le \delta} R_N(x,y)\, f(y)\, w(y)\, dy.
\]
Since we have established in \eqref{eq:AB_def} that $B_N(y)=1$, the expression reduces to
\[
Tf(x)= \frac{\gamma_{N+1}^2}{i}\, A_N(x)\int_{|y|\le \delta} f(y)\, w(y)\, dy\;+\;\int_{|y|\le \delta} R_N(x,y)\, f(y)\, w(y)\, dy.
\]
Now, by isolating each component we can decompose this as
\begin{equation}\label{eq:T_model}
T_{\mathrm{model}} f(x):= A_N(x)\int_{|y|\le \delta} f(y)\, w(y)\, dy,
\end{equation}
\begin{equation}\label{eq:T_smooth}
T_{\mathrm{smooth}} f(x):= \int_{|y|\le \delta} R_N(x,y)\, f(y)\, w(y)\, dy.
\end{equation}
Then adding \eqref{eq:T_model}, \eqref{eq:T_smooth}, we get the form
\begin{equation}\label{eq:T_decomposition}
Tf(x)= \frac{\gamma_{N+1}^2}{i}\, T_{\mathrm{model}} f(x)+ T_{\mathrm{smooth}} f(x).
\end{equation}
At this stage, it is possible to consider \eqref{eq:T_decomposition} for a reduction driven by its own internal structure. First of all, we can re-arrange \eqref{eq:T_decomposition} as the following:
\begin{equation}\label{eq:Lambda_def}
\Lambda f:= \int_{|y|\le \delta} f(y)\, w(y)\, dy.
\end{equation}
\begin{equation}\label{eq:MA_def}
(M_A c)(x):= A_N(x)\, c.
\end{equation}
\begin{equation}\label{eq:S_def}
(Sf)(x):= \int_{|y|\le \delta} R_N(x,y)\, f(y)\, w(y)\, dy.
\end{equation}
\begin{equation}\label{eq:T_compact}
T= \frac{\gamma_{N+1}^2}{i}\, M_A \circ \Lambda + S.
\end{equation}
From our construction \eqref{eq:AB_def}, $A_N$ is smooth on $[-\delta,\delta]$, hence bounded. Shrinking $\delta$ if necessary, we may assume
\begin{equation}\label{eq:AN_positive}
\inf_{|x|\le \delta} |A_N(x)| > 0.
\end{equation}
Consequently, $M_A$ is a bounded invertible operator on $L^p$. Since $R_N(x,y)$ is smooth and bounded on $(-\delta,\delta)^2$, the operator $S$ is bounded on $L^p(w)$ for any locally integrable weight $w$. From \eqref{eq:T_compact}, we write
\begin{equation}\label{eq:T_minus_S}
T - S = \frac{\gamma_{N+1}^2}{i}\, M_A \circ \Lambda.
\end{equation}
Composing on the left with the bounded inverse $M_A^{-1}$, we obtain the exact identity
\begin{equation}\label{eq:MA_inverse_identity}
M_A^{-1}(T - S) = \frac{\gamma_{N+1}^2}{i}\, \Lambda.
\end{equation}
Which pointwise in $x$ we have
\begin{equation}\label{eq:MA_inverse_Lambda_action}
\bigl(M_A^{-1}(T - S)f\bigr)(x)
= \frac{\gamma_{N+1}^2}{i}\, (\Lambda f).
\end{equation}
Now, by definition
\begin{equation}\label{eq:T_minus_S_action}
\bigl((T - S)f\bigr)(x)
= Tf(x) - Sf(x).
\end{equation}
So
\begin{equation}\label{eq:MA_inverse_action}
\bigl(M_A^{-1}(T - S)f\bigr)(x)
= \frac{1}{A_N(x)}\bigl(Tf(x) - Sf(x)\bigr).
\end{equation}
Substituting into \eqref{eq:MA_inverse_action}
\begin{equation}\label{eq:substituted_identity}
\frac{1}{A_N(x)}\bigl(Tf(x) - Sf(x)\bigr)
= \frac{\gamma_{N+1}^2}{i}\, (\Lambda f).
\end{equation}
Recalling \eqref{eq:S_def}, \eqref{eq:substituted_identity} becomes
\begin{equation}\label{eq:expanded_identity}
\frac{1}{A_N(x)}\left[Tf(x) - \int_{|y|\le \delta} R_N(x,y)\, f(y)\, w(y)\, dy \right] = \frac{\gamma_{N+1}^2}{i}\int_{|y|\le \delta} f(y)\, w(y)\, dy.
\end{equation}
Multiplying both sides of \eqref{eq:expanded_identity} by $\dfrac{i}{\gamma_{N+1}^2}$, we obtain
\begin{equation}\label{eq:scaled_identity}
\frac{i}{\gamma_{N+1}^2}\,\frac{1}{A_N(x)} \left[Tf(x) - \int_{|y|\le \delta} R_N(x,y)\, f(y)\, w(y)\, dy \right] = \int_{|y|\le \delta} f(y)\, w(y)\, dy.
\end{equation}
At this point, the right-hand side is a scalar depending only on $f$, and since the identity holds for every $x \in [-\delta,\delta]$, the left-hand side must be independent of $x$. Thus, for every $x$ with $|x|\le \delta$,
\begin{equation}\label{eq:final_identity_pointwise}
\int_{|y|\le \delta} f(y)\, w(y)\, dy = \frac{i}{\gamma_{N+1}^2}\, \frac{Tf(x) - Sf(x)}{A_N(x)}.
\end{equation}
We can conclude by defining the linear functional
\begin{equation}\label{eq:Lambda_def_final}
\Lambda f
:= \int_{|y|\le \delta} f(y)\, w(y)\, dy.
\end{equation}
Then, for all $|x|\le\delta$, one has the pointwise identity
\begin{equation}\label{eq:final_identity_pointwise}
\Lambda f = \frac{i}{\gamma_{N+1}^2}\, \frac{Tf(x) - Sf(x)}{A_N(x)} .
\end{equation}
From the reduction obtained in \eqref{eq:T_compact}–\eqref{eq:final_identity_pointwise}, boundedness of the weighted reconstruction operator is equivalent to boundedness of the scalar functional $\Lambda$ on $L^p(w)$. The problem therefore reduces to a purely measure-theoretic characterization of when $\Lambda$ defines a bounded linear functional on $L^p(w)$.

%%%% PROP ABOUT THE WEIGHTS (obs.1.ca.)

\begin{proposition}[Boundedness of the model functional]\label{prop:Lambda_boundedness}
Let $1 < p < \infty$, let $0 < \delta < \pi/4$, and let $w$ be a locally integrable weight on $(-\delta,\delta)$. Define the linear functional \eqref{eq:Lambda_def_final}
\[
\Lambda f := \int_{|y|\le \delta} f(y)\, w(y)\, dy.
\]
Then $\Lambda$ extends to a bounded linear functional on $L^p(w)$ if and only if
\begin{equation}\label{eq:Lambda_condition}
\int_{|y|\le \delta} w(y)^{-\frac{1}{p-1}}\, dy < \infty.
\end{equation}
In particular, if $w(y) \sim |y|^{\alpha}$ as $y \to 0$, then $\Lambda$ is bounded on $L^p(w)$ if and only if $\alpha < p - 1$.
\end{proposition}

\begin{proof}
Let $1 < p < \infty$ and let $w$ be a locally integrable weight on $(-\delta,\delta)$. Assume that
\[
\int_{|y|\le \delta} w(y)^{-\frac{1}{p-1}}\, dy < \infty.
\]
For any $f \in L^p(w)$, we apply Hölder's inequality with respect to the measure $w(y)\,dy$:
\[
|\Lambda f|
= \left| \int_{|y|\le \delta} f(y)\, w(y)\, dy \right|
\]
\[
= \left| \int_{|y|\le \delta} |f(y)|\, w(y)^{1/p}\, w(y)^{1-1/p}\, dy \right|
\]
\[
\le
\left( \int_{|y|\le \delta} |f(y)|^p\, w(y)\, dy \right)^{1/p}
\left( \int_{|y|\le \delta} w(y)^{1-p'}\, dy \right)^{1/p'} .
\]
Since $1 - p' = -\frac{1}{p-1}$, the second factor is finite by assumption, and hence
\[
|\Lambda f| \lesssim \|f\|_{L^p(w)}.
\]
Thus $\Lambda$ defines a bounded linear functional on $L^p(w)$. Conversely, suppose that $\Lambda$ is bounded on $L^p(w)$. Then $\Lambda \in (L^p(w))^{*}$, and by weighted duality one has
\[
(L^p(w))^{*} \simeq L^{p'}(w^{1-p'}),
\]
with pairing
\[
\langle f,h\rangle
= \int f(y)\, h(y)\, w(y)\, dy.
\]

Since $\Lambda f = \langle f,1\rangle$, boundedness of $\Lambda$ implies that the constant function $1$ belongs to $L^{p'}(w^{1-p'})$ on $(-\delta,\delta)$, that is,
\[
\int_{|y|\le \delta} w(y)^{1-p'}\, dy < \infty.
\]
Recalling again that $1 - p' = -\frac{1}{p-1}$, this leads to \eqref{eq:Lambda_condition}.
\end{proof}

Having established the boundedness criterion in Proposition~\ref{prop:Lambda_boundedness}, we now illustrate its applicability by considering weights with nontrivial local behavior near the mirror.

\paragraph{\textbf{\textit{Example A.}}}

Let $0<\delta\ll 1$ and define
\[
w(y)
:= |y|^{\alpha}\bigl(1+\lvert \sin(|y|^{-\gamma})\rvert\bigr)
\bigl(\log \tfrac{e}{|y|}\bigr)^{-\beta},
\qquad |y|\le \delta,
\]
where
\[
\alpha < p-1, \qquad \gamma>0, \qquad \beta \in \mathbb{R}.
\]
\textit{Verification.} By definition,
\[
w(y)^{-\frac{1}{p-1}}
=
|y|^{-\alpha/(p-1)}
\bigl(1+\lvert \sin(|y|^{-\gamma})\rvert\bigr)^{-\frac{1}{p-1}}
\bigl(\log \tfrac{e}{|y|}\bigr)^{\beta/(p-1)}.
\]
Now, for all $t$,
\[
1 \le 1+\lvert \sin t\rvert \le 2,
\]
hence
\[
2^{-\frac{1}{p-1}}
\le
\bigl(1+\lvert \sin t\rvert\bigr)^{-\frac{1}{p-1}}
\le 1.
\]

Therefore, there exist constants $c_1,c_2>0$ such that
\[
c_1\,|y|^{-\alpha/(p-1)}
\bigl(\log \tfrac{e}{|y|}\bigr)^{\beta/(p-1)}
\le
w(y)^{-\frac{1}{p-1}}
\le
c_2\,|y|^{-\alpha/(p-1)}
\bigl(\log \tfrac{e}{|y|}\bigr)^{\beta/(p-1)}.
\]

Thus,
\[
\int_{|y|\le \delta} w(y)^{-\frac{1}{p-1}}\,dy < \infty
\quad\Longleftrightarrow\quad
\int_{0}^{\delta}
y^{-\alpha/(p-1)}
\bigl(\log \tfrac{e}{y}\bigr)^{\beta/(p-1)}\,dy < \infty.
\]

Since $\alpha < p-1$, we have $-\alpha/(p-1) > -1$, and therefore
\[
\int_{0}^{\delta}
y^{-\alpha/(p-1)}
\bigl(\log \tfrac{e}{y}\bigr)^{\beta/(p-1)}\,dy < \infty
\qquad\text{for all } \beta \in \mathbb{R}.
\]

Hence
\[
w^{-\frac{1}{p-1}} \in L^1(0,\delta),
\]
and by Proposition~\ref{prop:Lambda_boundedness}, the functional
\[
\Lambda f := \int_{|y|\le \delta} f(y)\, w(y)\, dy
\]
is bounded on $L^p(w)$.

We have shown that, after localization and elimination of all smooth contributions, the weighted reconstruction operator reduces to a rank-one model functional. Boundedness on $L^p(w)$ is therefore governed entirely by the local integrability of $w^{-\frac{1}{p-1}}$ at the mirror, with no further dependence on oscillation, regularity, or higher-order degeneracy of the weight.

\end{document}